\documentclass{article}

\usepackage{amsmath,amsthm,amssymb,amsfonts}
\usepackage[vmargin=1in,hmargin=1.5in]{geometry}
\usepackage{graphicx}
\usepackage{tikz} \usetikzlibrary{cd} \usetikzlibrary{arrows}
\usepackage{pgfplots}
\usepackage{todonotes}
\usepackage[margin=40pt,font=small,labelfont=bf]{caption}
\usepackage{subcaption} 
\usepackage{hyperref}
\usepackage{cleveref}
\usepackage{rotating}
\usepackage{authblk}

\newtheorem{definition}{Definition}

\renewcommand{\max}{\textsf{Max}}
\newcommand{\city}{\textsf{City}}
\newcommand{\wikirandom}{\textsf{Wiki-random}}
\newcommand{\wrand}[1]{$\mathsf{WR_{#1}}$}
\newcommand{\random}{\textsf{Random}}
\newcommand{\rand}[1]{$\mathsf{R_{#1}}$}
\renewcommand{\P}{\mathcal{P}}

\pgfplotsset{compat=1.14}

\begin{document}

\title{\bf Topological Data Analysis on Simple English Wikipedia Articles}

\author[1]{Matthew Wright}
\affil[1]{St. Olaf College, \texttt{wright5@stolaf.edu}}

\author[2]{Xiaojun Zheng}
\affil[2]{Duke University, \texttt{xz264@duke.edu}}

\maketitle

\begin{abstract}
Single-parameter persistent homology, a key tool in topological data analysis, has been widely applied to data problems along with statistical techniques that quantify the significance of the results.
In contrast, statistical techniques for two-parameter persistence, while highly desirable  for real-world applications, have scarcely been considered.
We present three statistical approaches for comparing geometric data using two-parameter persistent homology; these approaches rely on the Hilbert function, matching distance, and barcodes obtained from two-parameter persistence modules computed from the point-cloud data.
Our statistical methods are broadly applicable for analysis of geometric data indexed by a real-valued parameter.
We apply these approaches to analyze high-dimensional point-cloud data obtained from Simple English Wikipedia articles.
In particular, we show how our methods can be utilized to distinguish certain subsets of the Wikipedia data and to compare with random data.
These results yield insights into the construction of null distributions and stability of our methods with respect to noisy data.
\end{abstract}

\section{Introduction and Motivation}

Topological Data Analysis (TDA) applies the mathematics of algebraic topology to study the shape of complex data with the goal of obtaining some insights about the data.
Persistent homology, a primary tool of TDA, is used to discern geometric and topological structure in high-dimensional datasets.
While single-parameter persistent homology has been widely used, multi-parameter variants of persistent homology are especially appealing, not only for their robustness in the presence of outliers, but also to analyze data naturally indexed by two or more parameters.
While computational techniques for two-parameter persistence have emerged in recent years \cite{LesnickWright}, there have yet been few practical applications and statistical techniques in this setting.

We present three statistical methods for two-parameter persistence and demonstrate the applicability of these methods to data arising from Simple English Wikipedia. 
Specifically, we consider large-scale hypothesis tests on Hilbert function values, the distribution of matching distance between two-parameter persistence modules, and the statistical significance of bar length.
We show that these methods are able to distinguish our Simple English Wikipedia from random point-cloud data.
We also consider the stability of these methods, specifically observing the distribution of matching distances as datasets become more dissimilar as well as the stability of Hilbert function values with respect to perturbations of the data.

\subsection{Project Data}

The data in this project was produced by applying a Word2Vec algorithm to the text of articles in Simple English Wikipedia, and was supplied to us by Shilad Sen \cite{Shilad}. 
The algorithm converted each of $120{,}526$ articles into a $200$-dimension vector, such that articles with similar content produce vectors that are close together (in the usual Euclidean metric).
The data also gives a \emph{popularity} score for each article, indicating how frequently the article is accessed in Simple English Wikipedia.
It is outside the scope of this article to consider the mechanics of producing the vectors and popularity scores; we simply treat the data as the input for our analyses. 
Abstractly, our data is a point cloud of $120{,}526$ points in $\mathbb{R}^{200}$, with a real-valued function on each point; our methods are applicable in this abstract setting.

\subsection{Software}

To obtain our results, we used the software RIVET, which computes and visualizes certain invariants of two-parameter persistence modules \cite{LesnickWright, rivet}.
RIVET takes geometric data as input and computes a two-parameter persistence module from the data. RIVET then outputs Hilbert function values of the module and barcodes along with user-requested linear slices of the module.
These barcodes allow us to compute matching distances between two such modules; these computations were performed using the \textsf{pyrivet} Python package \cite{rivetPython}, which depends on the Hera code for computing bottleneck distances \cite{hera, jea_hera}.
The \textsf{R} code for our statistical analyses is available at 
\url{https://github.com/Xiaojzheng/TDA-on-Simple-English-Wikipedia-Articles}.

\subsection{Outline}

The organization of the paper is as follows: In Section 2, we provide mathematical background related to two-parameter persistent homology. In Section 3, we describe the topological structures in the Simple English Wikipedia dataset and the three statistical measures that were performed on the dataset. We conclude this article with a discussion of results and directions for future research in Section 4.

\section{Mathematical Background}

Persistent homology is an algebraic method for identifying topological features such as connected components, holes, and voids in geometric data.
We briefly review persistent homology in the usual, single-parameter setting, and then generalize to the two-parameter setting---our focus.

\subsection{Single-Parameter Persistence}

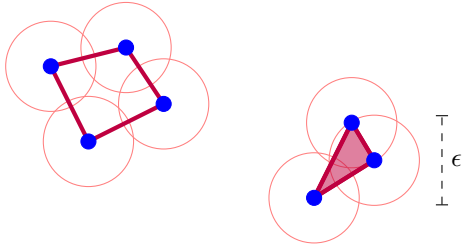
\begin{figure}[h]
  \centering
  \begin{tikzpicture}[scale=0.5]
    \foreach \Point in {(21,0),(15,1.5),(22,2),(14,3.5),(16,4),(22.6,1), (17,2.5)}{
      \draw[color=red!50] \Point circle (1.2);
    }
    \draw [purple, ultra thick] (14,3.5) -- (16,4);
    \draw [purple, ultra thick] (14,3.5) -- (15,1.5);
    \draw [purple, ultra thick] (22.6,1) -- (21,0);
    \draw [purple, ultra thick] (22,2) -- (22.6,1);
    \draw [purple, ultra thick] (21,0) -- (22,2);
    \draw [purple, ultra thick] (17,2.5) -- (16,4);
    \draw [purple, ultra thick] (15,1.5) -- (17,2.5);
    \draw[purple, ultra thick, fill=purple, fill opacity=0.5] (21,0) -- (22,2) -- (22.6,1) -- cycle;
    \foreach \Point in {(21,0),(15,1.5),(22,2),(14,3.5),(16,4), (22.6,1), (17,2.5)}{
      \filldraw[blue] \Point circle (0.2);
    }
    \draw [|-|, thin, dashed] (24.4,-0.2) -- node [right] {$\epsilon$} ++(0,2.4);

  \end{tikzpicture}

  \caption{Rips complex constructed from seven points at scale $\epsilon$. The circles of diameter $\epsilon$ are not part of the Rips complex but assist in visualising which pairs of points are within distance $\epsilon$.}
\label{complex}
\end{figure}

Single-parameter persistence studies the homology of a \textit{filtration}, which is a nested sequence of topological spaces. 
Often, these topological spaces are simplicial complexes constructed from a finite set of points, which we refer to as a \textit{point cloud}, in Euclidean space.
Given a point cloud $\P$ and scale parameter $\epsilon$, it is common to construct the \textit{Vietoris-Rips} (or simply \textit{Rips}) complex $R_\epsilon$, which consists of a $i$-simplex for every $(i+1)$ points of $\P$ whose pairwise distances are less than $\epsilon$, for each nonnegative integer $i$.
Figure \ref{complex} shows a Rips complex built from seven points.
 
If $\delta < \epsilon$, then $R_\delta \subseteq R_\epsilon$; thus the sequence of Rips complexes for increasing $\epsilon$ forms a filtration.
Assuming that $\P$ is finite, there are only finitely many distinct Rips complexes in the filtration. 
Thus, it suffices to consider a discrete filtration
\begin{equation}\label{filtration}
\P = R_0 \subset R_1 \subset R_2 \subset \cdots \subset R_{\mathrm{max}},
\end{equation}
where $R_{\mathrm{max}}$ is the complete simplex on $\P$.

Homology gives algebraic information about the topological features of a simplicial complex. We give a brief overview here; more details are found in a text such as \cite{CSGO} or \cite{hatcher2002algebraic}.
As is common in TDA, we compute homology over the two-element field $\mathbb{F}_2$.\footnote{Persistent homology can be defined with coefficients in any field, but $\mathbb{F}_2$ is often used in TDA for simplicity and ease of computation. The use of different fields yields different persistence diagrams due to the phenomena of torsion. The structure theorem for persistence modules requires coefficients from a field, rather than from a ring.}
Given a simplicial complex $X$, the set of all $i$-dimensional simplices in $X$ forms a basis for a vector space $C_i$ with coefficients in $\mathbb{F}_2$. The elements of $C_i$ are called \textit{$i$-chains}.
The \emph{boundary operator} $\partial_i \colon C_i \to C_{i-1}$ takes a $i$-simplex to the sum of its $(i-1)$-dimensional faces, extending by linearity to $i$-chains.
Let $B_i := \partial_{i+1}(C_{i+1}) \subseteq C_i$ be the subspace of boundaries, which are images of $\partial_{i+1}$.
Let $Z_i := \{ v \in C_i \mid \partial_i(v) = 0 \}$ be the subspace of cycles, which are $i$-chains with zero boundary.
Crucially, $B_i \subseteq Z_i$, since $\partial_i \circ \partial_{i+1} = 0$.
The \textit{homology vector space} is the quotient $H_i := Z_i/B_i$. 
Thus, $H_i$ is a vector space whose elements are equivalence classes of cycles that are not boundaries; two cycles are equivalent if they differ by a boundary.

The dimension of $H_i(X)$ is called the \textit{$i$th Betti number} of $X$, which is important for our purposes.
Intuitively, $H_0(X)$ is the number of connected components in $X$, $H_1(X)$ the number of holes, $H_2(X)$ the number of voids, and so on. If $X$ is the simplicial complex in Figure \ref{complex}, then $H_0(X) = 2$, $H_1(X) = 1$, and $H_i(X) = 0$ for $i \ge 2$, since the complex contains two connected components, one hole, and no higher homology.
In our work, we focus on degree-$0$ homology $H_0$ and degree-$1$ homology $H_1$.

For any index $k$, the subset relation $R_k \subset R_{k+1}$ in filtration \eqref{filtration} can be written as an inclusion map $R_k \hookrightarrow R_{k+1}$, which is a simplicial map. Accordingly, filtration \eqref{filtration} can be written as a sequence of simplicial maps:
\begin{equation}\label{simplicialMaps}
    R_0 \hookrightarrow R_1 \hookrightarrow R_2 \hookrightarrow \cdots \hookrightarrow R_{\mathrm{max}}.
\end{equation}
Importantly, homology is functorial, meaning that a simplicial map $X \to Y$ induces a linear map $H_i(X) \to H_i(Y)$. 
Functoriality allows us to take the homology of a filtration, obtaining a sequence of $\mathbb{F}_2$-vector spaces and linear maps.
For a fixed nonnegative integer $i$, the degree-$i$ homology of the sequence \eqref{simplicialMaps} is
\begin{equation}\label{persistenceModule}
  H_i(R_0) \to H_i(R_1) \to H_i(R_2) \to \cdots \to H_i(R_\mathrm{max}).
\end{equation}
Sequence \eqref{persistenceModule} may be regarded as infinite in both directions by prepending zero vector spaces and appending copies of $H_i(R_\mathrm{max})$, with identity maps between them.

\begin{figure}[ht]
  \centering
  \begin{tikzpicture}[scale=0.2]
    \node [rotate = 90] at (-13,3.8) {\small Filtration};
    \node [white] at (53,0) {.};
    
    \draw [thick] (-10,0) rectangle ++(7,7);
    \draw [thick] (1,0) rectangle ++(7,7);
    \draw [thick] (12,0) rectangle ++(7,7);
    \draw [thick] (23,0) rectangle ++(7,7);
    \draw [thick] (34,0) rectangle ++(7,7);
    \draw [thick] (45,0) rectangle ++(7,7);
  
    \filldraw [purple, opacity=0.4] (13,1) -- (13,2) -- (15,1);
    \filldraw [purple, opacity=0.4] (24,1) -- (24,2) -- (26,1);
    \filldraw [purple, opacity=0.4] (35,1) -- (35,2) -- (37,1);
    \filldraw [purple, opacity=0.4] (46,1) -- (46,2) -- (48,1);
    \filldraw [purple, opacity=0.4] (36,6) -- (37,1) -- (35,2);
    \filldraw [purple, opacity=0.4] (36,6) -- (40,4) -- (37,1);
    \filldraw [orange,opacity=0.5] (46,2) -- (48,1) -- (51,4) -- (47,6);
    
    \draw [purple, ultra thick] (2,1) -- (2,2);
    \draw [purple, ultra thick] (2,1) -- (4,1);
    \draw [purple, ultra thick] (13,1) -- (13,2);
    \draw [purple, ultra thick] (13,2) -- (15,1);
    \draw [purple, ultra thick] (24,1) -- (24,2);
    \draw [purple, ultra thick] (24,2) -- (26,1);
    \draw [purple, ultra thick] (35,1) -- (35,2);
    \draw [purple, ultra thick] (46,1) -- (46,2);
    \draw [purple, ultra thick] (46,2) -- (48,1);
    \draw [purple, ultra thick] (13,1) -- (15,1);
    \draw [purple, ultra thick] (24,1) -- (26,1);
    \draw [purple, ultra thick] (35,1) -- (37,1);
    \draw [purple, ultra thick] (46,1) -- (48,1);
    \draw [purple, ultra thick] (24,2) -- (25,6);
    \draw [purple, ultra thick] (26,1) -- (29,4);
    \draw [purple, ultra thick] (25,6) -- (29,4);
    \draw [purple, ultra thick] (36,6) -- (40,4);
    \draw [purple, ultra thick] (47,6) -- (51,4);
    \draw [purple, ultra thick] (37,1) -- (40,4);
    \draw [purple, ultra thick] (35,2) -- (36,6);
    \draw [purple, ultra thick] (37,1) -- (35,2);
    \draw [purple, ultra thick] (36,6) -- (37,1);
    \draw [purple, ultra thick] (48,1) -- (47,6);
    \draw [purple, ultra thick] (46,2) -- (47,6);
    \draw [purple, ultra thick] (46,2) -- (51,4);
    \draw [purple, ultra thick] (48,1) -- (51,4);
    
    \node at (-6.5,8.8) {$R_1$};
    \node at (4.5,8.8) {$R_2$};
    \node at (15.5,8.8) {$R_3$};
    \node at (26.5,8.8) {$R_4$};
    \node at (37.5,8.8) {$R_5$};
    \node  at (48.5,8.8) {$R_6$};
    
    \draw [thick, right hook-latex, red] (-4.5, 8.7) -- ++(6.8,0);
    \draw [thick, right hook-latex, red] (6.5, 8.7) -- ++(6.8,0);
    \draw [thick, right hook-latex, red] (17.5, 8.7) -- ++(6.8,0);
    \draw [thick, right hook-latex, red] (28.5, 8.7) -- ++(6.8,0);
    \draw [thick, right hook-latex, red] (39.5, 8.7) -- ++(6.8,0);
    
    \foreach \Point in {(-9,1),(-9,2),(-7,1),(-8,6),(-4,4),(2,1),(2,2),(4,1),(3,6),(7,4),(13,1),(13,2),(15,1),(14,6),(18,4),(24,1),(24,2),(26,1),(25,6),(29,4),(35,1),(35,2),(37,1),(36,6),(40,4),(46,1),(46,2),(48,1),(47,6),(51,4)}{
        \filldraw[blue] \Point circle (0.32);
    }
  \end{tikzpicture}
  \begin{tikzpicture}[scale=1.11]
    \draw [line width=0.1cm,orange] (0.01,-4.2) -- (1.6,-4.2);
    \draw [line width=0.1cm,orange] (0.01,-4.5) -- (2,-4.5);
    \draw [line width=0.1cm,orange] (0.01,-4.8) -- (5.5,-4.8);
    \draw [line width=0.1cm,orange] (0.01,-5.1) -- (5.8,-5.1);
    \draw [line width=0.1cm,orange] (0.01,-5.4) -- (11,-5.4);
    \draw [line width=0.1cm,blue] (5.9,-5.7) -- (8.1,-5.7);
    \node [rotate = 90] at (-0.4,-5) {\small Barcode};

    \draw [dashed] (0.7,-3.8) --(0.7,-6.1);
    \draw [dashed] (2.7,-3.8) --(2.7,-6.1);
    \draw [dashed] (4.7,-3.8) --(4.7,-6.1);
    \draw [dashed] (6.7,-3.8) --(6.7,-6.1);
    \draw [dashed] (8.7,-3.8) -- (8.7,-6.1);
    \draw [dashed] (10.7,-3.8) -- (10.7,-6.1);

    \draw [->,>=stealth] (0,-6.1) -- (11.3,-6.1);
    \draw [line width=0.02cm] (0,-6.1) -- (0,-3.8);
    \node [scale = 0.03cm] at (11.5,-6.1) {$\epsilon$};
  \end{tikzpicture}

\caption{Six complexes are shown from a Rips filtration built from a five points. In the corresponding barcode, orange bars denote $H_0$ homology (indicating connected components) and blue bars denote $H_1$ homology (indicating holes).}
\label{barcode}
\end{figure}
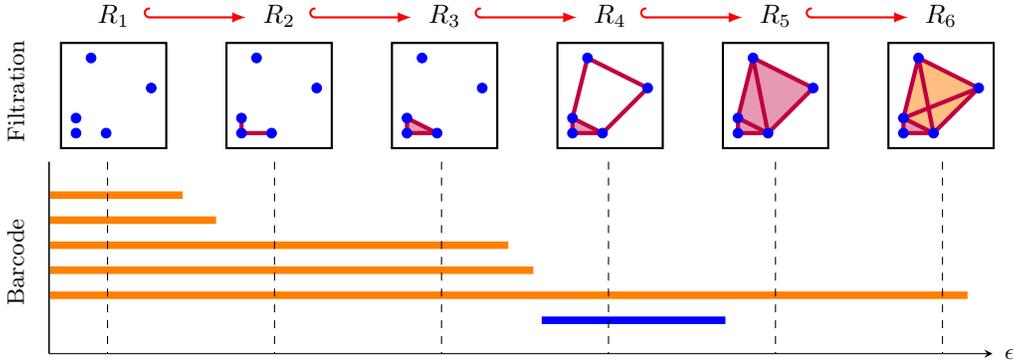

The homology of a filtration, as in sequence \eqref{persistenceModule}, is a \textit{persistence module}.
As an algebraic structure, a persistence module is an algebraic module over the polynomial ring $\mathbb{F}_2[x]$, where the action of $x$ shifts module elements forward in the sequence, taking elements of $H_i(R_k)$ to elements of $H_i(R_{k+1})$ via the linear maps in sequence \eqref{persistenceModule}.
By the structure theorem for finitely-generated modules over principal ideal domains, a persistence module decomposes as a sum of interval modules. An interval module consists of copies of the field $\mathbb{F}_2$ with identity maps between them at some consecutive indexes and zero vector spaces at all other indexes:
\[ \cdots \to 0 \to \mathbb{F}_2 \to \mathbb{F}_2 \to \cdots \to \mathbb{F}_2 \to \mathbb{F}_2 \to 0 \to \cdots. \]
Each interval module encodes the birth and death of one $i$-dimensional topological feature of the filtration, with the field $\mathbb{F}_2$ appearing at exactly the filtration indexes at which the feature is present.
Taken together, the interval modules give the \textit{barcode} of the filtration: a multiset of intervals, each of which indicates the lifespan of one $i$-dimensional feature in the filtration.
Figure \ref{barcode} shows an example of a filtration and its barcode, adapted from \cite{HanOkoYadZhe}.
For more details on persistent homology, see \cite{ghristbarcodes} or \cite{CSGO}.

A barcode may also be visualized by a \textit{persistence diagram}, in which each bar is plotted as an ordered pair.
Thus, a persistence diagram consists of a multiset of points $(b,d)$ with $b \le d$. Here, $b$ denotes the birth and $d$ the death of a topological feature; the quantity $d-b$ is the lifespan of the feature.

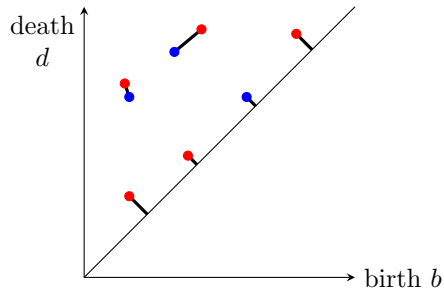
\begin{figure}[ht]
  \centering

  \begin{tikzpicture}[scale=0.6]
    \draw [<->,>=stealth] (7,1) node[right]{birth $b$} -- (1,1) -- (1,7) node[below left, align=center]{death\\$d$};
    \draw (1,1) -- (7,7);

    \draw [very thick,black] (3,6) -- (3.6,6.5);
    \draw [very thick,black] (2,5) -- (1.9,5.3);
    \draw [very thick,black] (2.4,2.4) -- (2,2.8);
    \draw [very thick,black] (3.3,3.7) -- (3.5,3.5);
    \draw [very thick,black] (5.7,6.4) -- (6.05,6.05);
    \draw [very thick,black] (4.6,5) -- (4.8, 4.8);

    \foreach \Point in {(3,6),(2,5),(4.6,5)}{
        \filldraw [blue] \Point circle (0.1);
    }
    \foreach \Point in {(3.6,6.5),(1.9,5.3),(2,2.8),(3.3,3.7),(5.7,6.4)}{
        \filldraw [red] \Point circle (0.1);
    }
  \end{tikzpicture}

  \caption{A matching between persistence diagrams $\mathcal{D}_1$ (plotted in blue) and $\mathcal{D}_2$ (plotted in red). Some points in each diagram are matched to the diagonal. If this matching minimizes the $L_\infty$ distance between matched points, then its greatest $L_\infty$ distance between matched points is the bottleneck distance between $\mathcal{D}_1$ and $\mathcal{D}_2$.}
\label{matching}
\end{figure}

We compare two persistence diagrams via the bottleneck distance, which can be thought of as the cost of transforming one diagram into the other.
To define this distance, we must introduce the concept of a matching between persistence diagrams: this is a sort of bijection between diagrams, but we allow points in a diagram to be paired with points on the diagonal. 
To explain this, let $\mathcal{L}$ be the diagonal line consisting of all points $(b,b)$.
A \textit{matching} between persistence diagrams $\mathcal{D}_1$ and $\mathcal{D}_2$ is a bijection between $\mathcal{D}_1 \cup \mathcal{L}$ and $\mathcal{D}_2 \cup \mathcal{L}$.
\Cref{matching} illustrates a matching between two persistence diagrams. 

Given a matching, we compute the maximum $L_\infty$ distance between all pairs of matched points. 
The smallest such distance, computed over all possible matchings of diagrams $\mathcal{D}_1$ and $\mathcal{D}_2$, gives the bottleneck distance between $\mathcal{D}_1$ and $\mathcal{D}_2$.
Intuitively, the bottleneck distance is the largest adjustment of birth or death coordinates required to transform diagram $\mathcal{D}_1$ into $\mathcal{D}_2$ (or vice-versa). Formally, the bottleneck distance is defined as follows.

\begin{definition} 
The \emph{bottleneck distance} between persistence diagrams $\mathcal{D}_1$ and $\mathcal{D}_2$ is
\[ d_B(\mathcal{D}_1, \mathcal{D}_2) = \inf_\eta \sup_x || x - \eta(x) ||_\infty, \]
where the infimum is taken over all matchings $\eta$ of $\mathcal{D}_1$ and $\mathcal{D}_2$, and the supremum is taken over all pairs of matched points $x$ and $\eta(x)$.
\end{definition}

We regard the bottleneck distance between $\mathcal{D}_1$ and $\mathcal{D}_2$ as a topological notion of distance between the point cloud datasets from which $\mathcal{D}_1$ and $\mathcal{D}_2$ are computed.

\subsection{Two-parameter persistence}\label{twoparam}

A \textit{bifiltration} is a two-dimensional generalization of a filtration: a set of simplicial complexes, each indexed by two parameters, with commuting inclusion maps in the direction of increase of each parameter. 
Bifiltrations are useful for data simultaneously indexed by two parameters, as described in \cite{CarlssonZomorodian}; we explain one common construction here.

Given a point cloud $\P \subset \mathbb{R}^d$ and a function $\gamma : \P \to \mathbb{R}$, we define the \textit{function-Rips} bifiltration, which we will use in sequel. 
For a pair of parameters $\alpha$ and $\epsilon$, let $R_{\alpha,\epsilon}$ be the Rips complex with scale parameter $\epsilon$ constructed from $\gamma^{-1}((-\infty,\alpha])$, the subset of points $p \in \P$ such that $\gamma(p) \le \alpha$. 
Thus, we obtain a bi-indexed set of Rips complexes such that $\delta < \epsilon$ implies $R_{\alpha,\delta} \subseteq R_{\alpha,\epsilon}$, and also $\beta < \alpha$ implies $R_{\beta,\epsilon} \subseteq R_{\alpha,\epsilon}$.
Assuming again that $\P$ is finite, simplices appear at a discrete set of indexes, so it suffices to consider a discrete bifiltration which we can re-index with integer indexes, as in the following commutative diagram.

\begin{center}
  \begin{tikzcd}
    \vdots & \vdots & \vdots \\
    R_{0,2} \arrow[r, hook] \arrow[u, hook] & R_{1,2} \arrow[r, hook] \arrow[u, hook] & R_{2,2} \arrow[r, hook] \arrow[u, hook] & \cdots \\
    R_{0,1} \arrow[r, hook] \arrow[u, hook] & R_{1,1} \arrow[r, hook] \arrow[u, hook] & R_{2,1} \arrow[r, hook] \arrow[u, hook] & \cdots \\
    R_{0,0} \arrow[r, hook] \arrow[u, hook] & R_{1,0} \arrow[r, hook] \arrow[u, hook] & R_{2,0} \arrow[r, hook] \arrow[u, hook] & \cdots 
  \end{tikzcd}
\end{center}

Though finite $\P$ implies only finitely many distinct complexes in the diagram above, it is often convenient to regard the diagram above as infinite. Following the arrows to the right or up eventually leads to a constant sequence of complexes and identity maps. We may further regard the diagram as infinite in all directions by supplying empty complexes at negative indexes.

The homology of a bifiltration is a \textit{two-parameter persistence module}: a bi-indexed set of vector spaces, with commuting linear maps between them as in the following diagram. Again, this diagram is constant for sufficiently large index values, and may be extended to negative indexes with zero vector spaces.

\begin{center}
  \begin{tikzcd}
    \vdots & \vdots & \vdots \\
    H_i(R_{0,2}) \arrow[r] \arrow[u] & H_i(R_{1,2}) \arrow[r] \arrow[u] & H_i(R_{2,2}) \arrow[r] \arrow[u] & \cdots \\
    H_i(R_{0,1}) \arrow[r] \arrow[u] & H_i(R_{1,1}) \arrow[r] \arrow[u] & H_i(R_{2,1}) \arrow[r] \arrow[u] & \cdots \\
    H_i(R_{0,0}) \arrow[r] \arrow[u] & H_i(R_{1,0}) \arrow[r] \arrow[u] & H_i(R_{2,0}) \arrow[r] \arrow[u] & \cdots
  \end{tikzcd}
\end{center}

The algebraic structure of two-parameter persistence modules is quite complicated.
As such, there is no analog of a barcode for two-parameter persistence modules; for more details see \cite{lesnick}. Instead, we obtain a barcode along each ``slice'' with nonnegative slope through the two-parameter persistence module, which we now explain.

Defining barcodes along slices of a two-parameter persistence module $\mathcal{M}$ requires returning to the continuous perspective.
We work in the $\alpha\epsilon$-plane, where $\alpha \in \mathbb{R}$ is the function value parameter and $\epsilon \ge 0$ is the scale parameter.
For any pair of parameters $(\alpha, \epsilon)$, there is a Rips complex $R_{\alpha, \epsilon}$ with degree-$i$ homology $H_i(R_{\alpha, \epsilon})$. 
If $(\beta, \delta)$ is such that $\beta \le \alpha$ and $\delta \le \epsilon$, then there is a linear map $H_i(R_{\beta, \delta}) \to H_i(R_{\alpha, \epsilon})$.

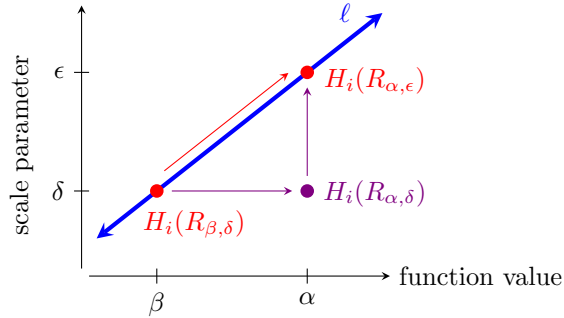
\begin{figure}[h]
  \centering
  \begin{tikzpicture}
    \draw [->,>=stealth] (0.1,0) -- ++(4,0) node[right] {function value};
    \draw (1,-0.1) node[below]{$\beta$} -- (1,0.1);
    \draw (3,-0.1) node[below]{$\alpha$} -- (3,0.1);
    
    \draw [->,>=stealth] (0,0.1) -- ++(0,3.5);
    \draw (0.1,1.132) -- (-0.1,1.132) node[left] {$\delta$};
    \draw (0.1,2.71) -- (-0.1,2.71) node[left] {$\epsilon$};
    \node [rotate=90] at (-0.8,1.75) {scale parameter};
    
    \draw[<->, blue, ultra thick,>=stealth] (0.2,0.5) -- ++(3.8,3);
    \node [blue] at (3.5, 3.5) {$\ell$};
    
    \fill[red] (1,1.132) circle (0.09);
    \fill[red] (3,2.71) circle (0.09);
    \fill[violet] (3,1.132) circle (0.09);
    
    \node[red] at (1.5,0.7) {$H_i(R_{\beta, \delta})$};
    \node[violet] at (3.9,1.1) {$H_i(R_{\alpha, \delta})$};
    \node[red] at (3.9,2.6) {$H_i(R_{\alpha, \epsilon})$};
    
    \draw[->,>=stealth,violet] (1.2,1.132) -- ++ (1.6,0);
    \draw[->,>=stealth,violet] (3,1.332) -- (3,2.51);
    \draw[->,>=stealth,red] (1.1,1.4) -- (2.75,2.703);
  \end{tikzpicture}
  \caption{Given a two-parameter persistence module $\mathcal{M}$ and a line $\ell$ of nonnegative slope, we define a one-parameter persistence module $\mathcal{M}_\ell$ consisting of homology vector spaces from $\mathcal{M}$ at points in $\ell$, with their induced linear maps.}
  \label{sliceModule}
\end{figure}

Let $\ell$ be a line of nonnegative slope in the $\alpha\epsilon$-plane.\footnote{The line $\ell$ must be parameterized in order to make sense of persistence diagrams along $\ell$. As described in \cite{LesnickWright}, we choose the parameterization $\gamma_\ell : \mathbb{R} \to \ell$ to be the unique order-preserving isometry such that $\gamma_\ell(0)$ is the point where $\ell$ intersects the non-negative portions of the coordinate axes.}
We define a one-parameter persistence module $\mathcal{M}_\ell$ by taking the ``slice'' of $\mathcal{M}$ along $\ell$. Specifically, to each point $(\alpha, \epsilon)$ in line $\ell$ we assign the corresponding homology vector space $H_i(R_{\alpha, \epsilon})$, with linear maps induced by those in $\mathcal{M}$, as illustrated in \Cref{sliceModule}. 
Thus,  $\mathcal{M}_\ell$ is a one-parameter family of homology vector spaces and linear maps. 
As a one-parameter persistence module, $\mathcal{M}_\ell$ has a barcode, or equivalently, a persistence diagram.
More details on slices of two-parameter persistence modules can be found in the friendly introduction \cite{HanOkoYadZhe} or the comprehensive \cite{LesnickWright}.

Next, we define a distance between two-parameter persistence modules based on the bottleneck distance between persistence diagrams along all possible lines of nonnegative slope in the $\alpha\epsilon$-plane.

\begin{definition}\label{matchingDist}
The \emph{matching distance}, $d_{M}$, between two-parameter persistence modules $\mathcal{M}$ and $\mathcal{N}$ is the supremum of the bottleneck distances between the persistence diagrams along corresponding lines of non-negative slope in the two modules. 
Precisely,
\[ d_{M} = \sup_{\ell} \{ d_B(\mathcal{D}(\mathcal{M}_\ell), \mathcal{D}(\mathcal{N}_\ell)) \cdot w(\mathrm{slope}(\ell)) \}, \]
where the supremum is over all lines of nonnegative slope and $w(m) = \frac{1}{\sqrt{1+q^2}}$, where $q = \mathrm{max}\left(m, \frac{1}{m}\right)$.

\end{definition} 

We note that in Definition \ref{matchingDist}, a weight $w$ is assigned to each line, which depends on the slope $\ell$. A line with slope one is assigned the largest weight, and the weight approaches zero as the slope approaches zero or infinity. This gives greatest weight to slopes at which the persistence diagrams are most stable.\footnote{Furthermore, the weight is chosen such that if the interleaving distance between persistence modules $\mathcal{M}$ and $\mathcal{N}$ is $1$, then the weighted bottleneck distance is at most $1$. For details, see \cite{Landi2018} or \cite{lesnick}.}

We now introduce two important functions on two-parameter persistence modules.
First, the \textit{Hilbert function} gives the dimension of the homology vector space $H_i(R_{\alpha, \epsilon})$ at every pair of parameter values $(\alpha,\epsilon)$. Thus, the Hilbert function provides a useful summary of the structure of the bifiltration. For $i=0$, the Hilbert function gives the number of connected components in $R_{\alpha,\epsilon}$. For $i=1$, the Hilbert function gives the number of holes (or $1$-cycles) in $R_{\alpha,\epsilon}$. 

Second, the \textit{bigraded Betti numbers} count the number of topological features (i.e., components for $i=0$ and holes for $i=1$) that are born or die at each pair of parameter values $(\alpha,\epsilon)$. The zeroth bigraded Betti number, denoted $\xi_0(\alpha, \epsilon)$, counts the number of topological features  that are born at parameter pair $(\alpha,\epsilon)$. In other words, $\xi_0(\alpha,\epsilon)$ is the number of homology classes that exist in $H_i(R_{\alpha,\epsilon})$ but not in $H_i(R_{\beta,\delta})$ if either $\beta < \alpha$ or $\delta < \epsilon$. 
Similarly, the first bigraded Betti number, denoted $\xi_1(\alpha,\epsilon)$, counts the number of topological features that die at parameter pair $(\alpha,\epsilon)$ --- homology classes that exist in $H_i(R_{\beta,\delta})$ if $\beta < \alpha$ and $\delta < \epsilon$, but not in $H_i(R_{\alpha,\epsilon})$.
Algebraically, $\xi_0$ counts the generators and $\xi_1$ counts the relations of the persistent homology module at each pair of parameter values. For more details, see \cite{LesnickWright}.

\begin{figure}[h]
    \centering
    \includegraphics[scale=0.6]{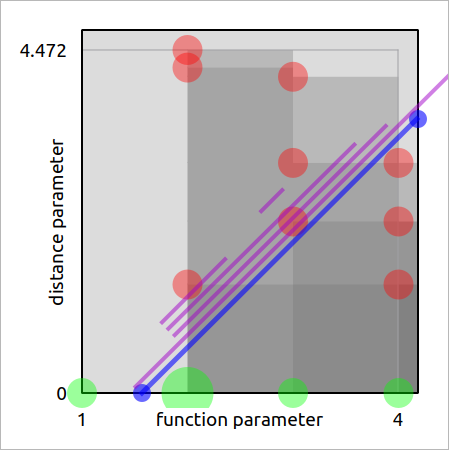}
    \caption{An example RIVET plot. 
    }
    \label{rivet_plot}
\end{figure}

The RIVET software computes and visualizes both the Hilbert function values and the bigraded Betti numbers for a two-parameter persistence module. A sample RIVET visualization is shown in \Cref{rivet_plot}.
The plot shows the two-parameter plane, with the function value on the horizontal axis, and the distance scale on the vertical axis. 
The Hilbert function values are depicted by grayscale shading, with darker shades of gray indicating larger values of the Hilbert function.
The bigraded Betti numbers are shown as green and red dots: green dots for $\xi_0$ and red dots for $\xi_1$. The area of the dots is proportional to the value of $\xi_0$ or $\xi_1$ --- the number of topological features that are born or die --- at each point in the two-parameter space.
RIVET also displays a barcode along a user-selected slice of the two-parameter persistence module: the slice line is shown in blue, and the barcode in purple.

When computing with RIVET, the user must specify the level of discretization, colloquially referred to as the number of \emph{bins}, to be applied to each coordinate axis. A small number of bins results in a faster, coarser calculation, while a large number of bins yields finer detail but with considerably longer runtime. The choice of bins determines the resolution at which we view the module: if we select $m$ horizontal bins and $n$ vertical bins, then we obtain Hilbert function values and bigraded Betti numbers at $mn$ points in the two-parameter plane.
For most of our calculations, we choose $20$ bins along each coordinate axis.

A similar discretization issue arises when approximating a matching distance between two-parameter persistence modules using \textsf{pyrivet}. The algorithm approximates the matching distance by computing barcodes along a finite set of lines; each line is determined by its angle from the horizontal and its offset from the origin. The user must specify how many angles and offset values to use in the calculation; more values generally result in a better approximation, but with the cost of additional computation time. After some preliminary tests, we settled on $20$ angles and $20$ offset values to obtain accurate results in reasonable computation time.

\section{Data Analysis}

\subsection{Sampling}

Our Wikipedia dataset consists of $120{,}526$ vectors in $\mathbb{R}^{200}$; each vector is obtained by applying a Word2Vec algorithm to a Simple English Wikipedia article. We regard the distance between vectors as indicative of the semantic relatedness of the corresponding articles \cite{Shilad}. 
Each article vector was accompanied by a popularity value, with larger values indicating more popular articles.

We selected three subsets of the Wikipedia dataset, which we refer to as \max, \city, and \wikirandom. 
The \max\ dataset consists of the $2000$ vectors with the largest popularity values.
The \city\ dataset consists of $2000$ vectors randomly chosen from the vectors whose article title is a city name. 
The \wikirandom\ dataset contains $2000$ randomly chosen vectors from the Wikipedia dataset.
Part of our analysis required many independent random subsets of the Wikipedia data; we refer to these as \wrand{1}, \wrand{2}, and so on, each containing  $2000$ randomly-chosen vectors from the Wikipedia dataset.
Within each dataset, we sorted the vectors in order of decreasing popularity; the index of each vector in the sorted list then gives a ``popularity rank'' which serves as our real-valued function on the dataset. Thus, each function-Rips bifiltration was constructed using the popularity rank indexes---integers from 1 to 2000 for the vectors in each dataset. Importantly, lower popularity rank indexes indicate \emph{more popular} articles, which appear first in the bifiltration.

Furthermore, we created many random datasets consisting of $2000$ vectors in $\mathbb{R}^{200}$.
We observed that the distribution of components of the Wikipedia vectors is very close to a normal distribution; thus we sampled from this distribution to obtain components for our random vectors. Furthermore, we randomly assigned ``popularity rank'' values (integers from 1 to 2000) to the vectors in each random dataset.
A single such dataset is denoted \random; a sequence of independently-generated random data sets is denoted \rand{1}, \rand{2}, etc.
We now compare the topological features between subsets of the Wikipedia data and contrast these with sets of random vectors.

\subsection{Preliminary Observations}

\begin{sidewaysfigure}
    \centering
    	\includegraphics[scale=0.525]{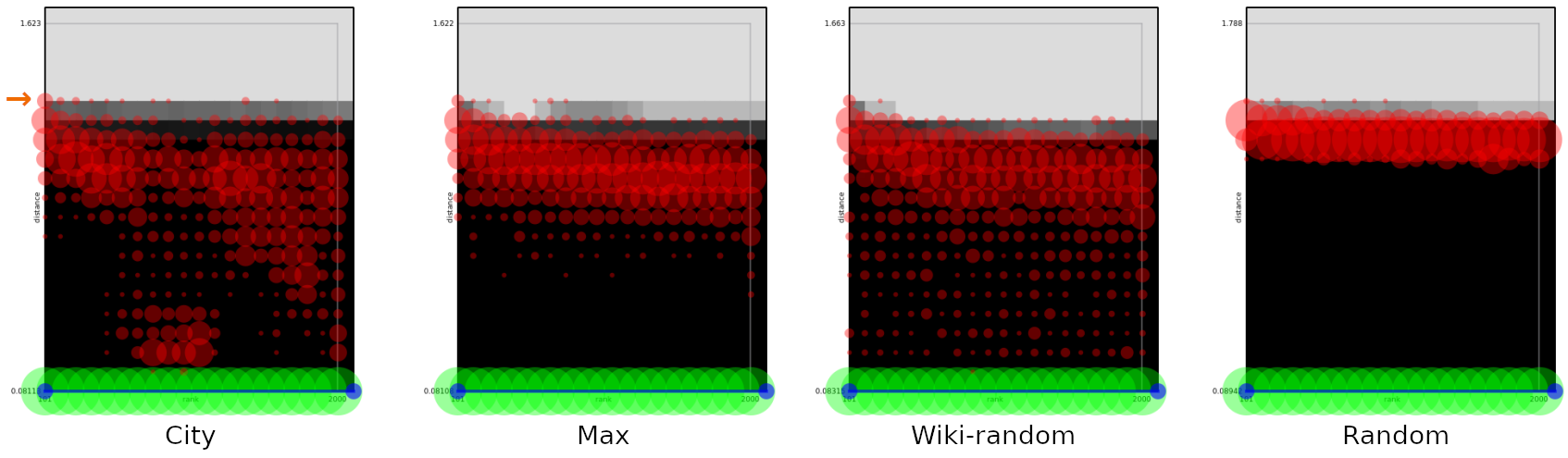}
    	\caption{RIVET plots of $H_0$ persistence modules of four datasets.}
    \label{wiki_H0_all}
    \vspace{1cm}
    	\includegraphics[scale=0.67]{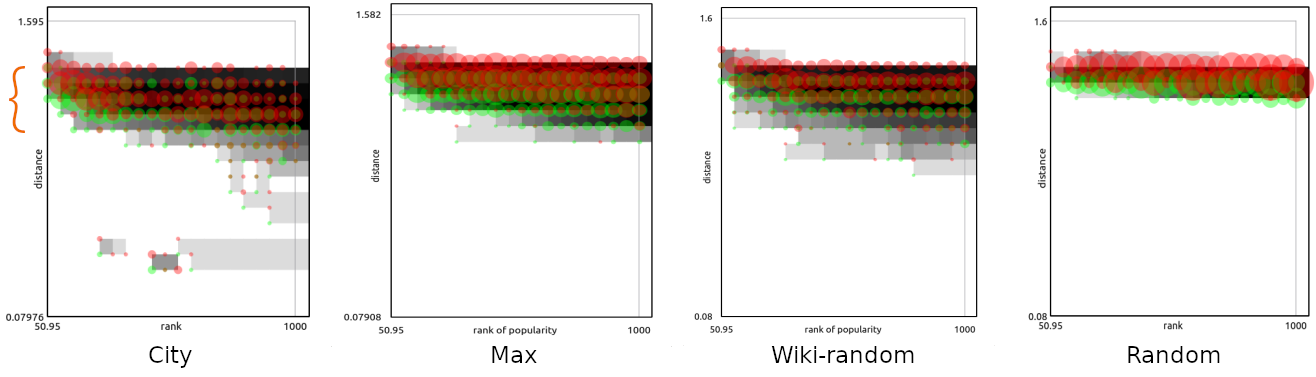}
    \caption{RIVET plots of $H_1$ persistence modules of four datasets.}
    \label{wiki_H1_all}
\end{sidewaysfigure}

\Cref{wiki_H0_all} shows the Hilbert function and bigraded Betti number visualization produced by RIVET for the degree-$0$ homology ($H_0$) of the four data sets; analogous plots for degree-$1$ homology ($H_1$) appear in \Cref{wiki_H1_all}.

The $H_0$ plots show that the Hilbert function has large values when the scale parameter is between zero and a certain threshold, indicated roughly by the orange arrow at left in \Cref{wiki_H0_all}, above which the Hilbert function value is one. That is, nontrivial degree-0 homology only exists when the scale parameter is less than this threshold; above the threshold the entire data set is connected into one large component. 

Furthermore, in the \city\ and \wikirandom\ plots, we notice red dots at small values on the vertical axis, indicating that some data points are quite close together. This is especially noticeable in the \city\ data: there exist certain clusters of cities whose articles are close together in the semantic space.
Interestingly, in the \random\ dataset, the red dots occur in a narrow band, indicating little variation in the distance from a point to its neighbors.
The distance scale at which components in \random\ connect is larger than the distance scale at which components connect in the other data sets. Furthermore, all of the connections in \random\ occur within a narrow range of distance values.

In the $H_1$ plots, the Hilbert function value is zero for small values of the scale parameter, then large in a certain interval of scale values, before becoming zero again. We again see that nontrivial degree-1 homology exists mostly in a narrow range of parameter values, indicated roughly by the orange bracket at left in \Cref{wiki_H1_all}, which is slightly above the range of distance values in which degree-0 homology exists.
This mirrors a phenomenon seen elsewhere in random topology: each degree of homology dominates the others in a certain interval of scale \cite{BobrowskiKahle}.
Indeed, this phenomenon is most pronounced for our \random\ dataset, for which the nontrival degree-1 homology occurs in a range narrower than the other datasets.

Having observed structure in the RIVET plots, we turn to statistical analysis of the topological difference between the Wikipedia data sets and the random data sets.
We assume that the Wikipedia datasets contain information not present in the random datasets, which form our null distribution. We will look for statistical significance to tell us that a dataset is not likely to be random.

\subsection{Statistical Methods}

\subsubsection{Large Scale Hypothesis Testing}\label{largeScaleTesting}

The first statistical technique that we used is \emph{large scale hypothesis testing}. 
A hypothesis test evaluates the probability of observing the given data based on some null hypothesis about the distribution from which the data has been sampled. If this probability is below a specified threshold, then the null hypothesis is rejected.
For more background on hypothesis testing, see \cite{Efron} or \cite{Wasserman}. When thousands of hypothesis tests are to be conducted, large scale hypothesis testing can be employed to decide which null hypotheses to reject. 
For example, large scale testing is used to compare gene expression levels between cancer patients and a control group; one hypothesis test is performed for each of thousands of genes under consideration \cite{Efron}. 
We use large scale testing to compare Hilbert function values at corresponding pixels across multiple RIVET plots. 
We assume that pixels at the same coordinate in different plots convey information about the same topological property in different datasets, so the pixel values at the same coordinate across various plots form a single hypothesis test. Our goal is to identify which pixels differ significantly between the Wikipedia datasets and random datasets.

The RIVET computation outputs Hilbert function values, visualized in grayscale in the plots in \Cref{wiki_H0_all} and \Cref{wiki_H1_all}.
Each of these values gives the number of connected components (for $H_0$ homology) or holes (for $H_1$) for a certain popularity threshold and distance scale. Since we used $20$ bins for each parameter, we obtained $400$ values for each dataset in each homology degree.
This results in $400$ hypothesis tests for $H_0$ homology, and another $400$ hypothesis tests for $H_1$ homology.

We applied large scale hypothesis testing to compare our Wiki-random and random datasets.
Specifically, we used 15 samples of Wiki-random data \wrand{1}, \wrand{2}, ..., \wrand{15}, and 15 samples of random data \rand{1}, \rand{2}, ..., \rand{15}.

\begin{figure}[h]
    \centering
    \includegraphics[scale=0.5]{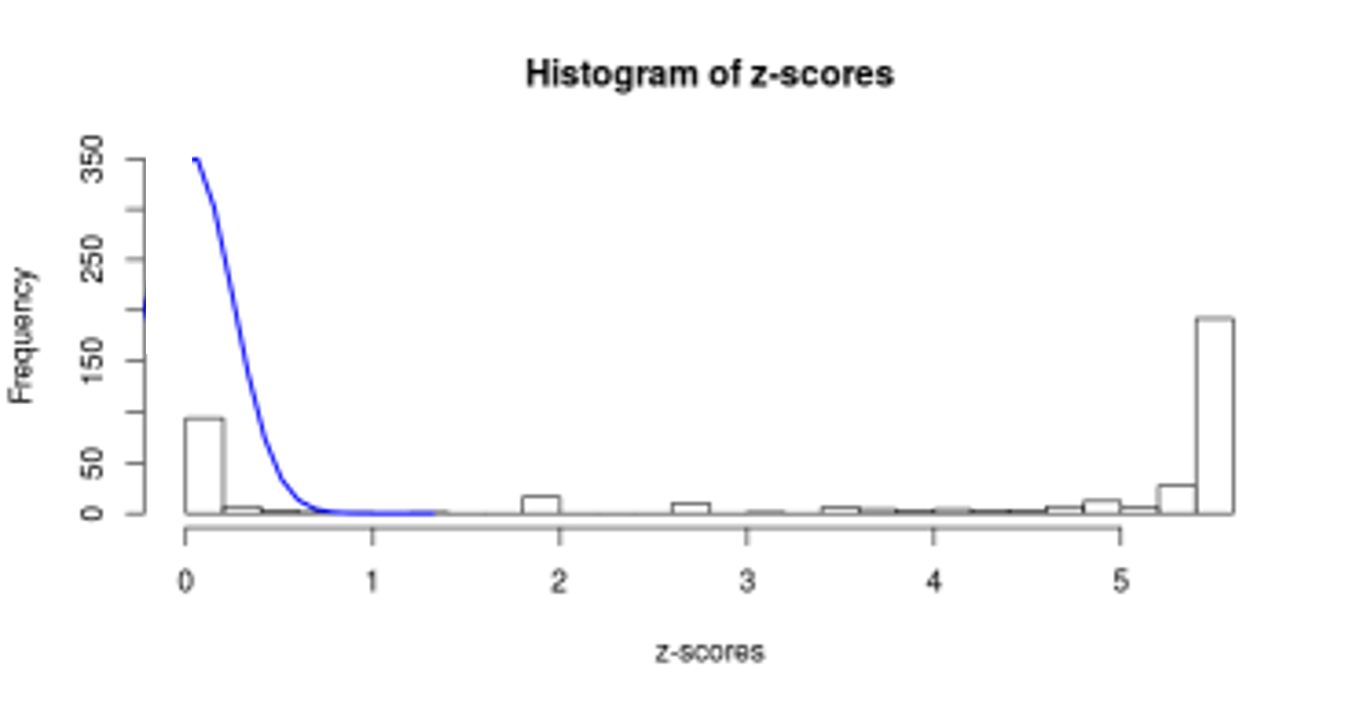}
    \caption{The null distribution of z-scores for the large scale hypothesis test is plotted in blue, and the histogram shows the actual z-scores when comparing the Wiki-random and random datasets.}
    \label{Hist_LST}
\end{figure}

We first examined the Hilbert function values for the $H_0$ homology computed from our datasets.
Using the random datasets, we applied cross-validation to create a null distribution. We ran $500$ experiments; in each, we randomly partitioned \rand{1}, \rand{2}, ..., \rand{15} into two subsets of size $7$ and $8$.
For each of the $400$ pixels, we compared the two subsets by a two-sample t-test. This resulted in $500$ t-statistics for each pixel. Finally, we took the average of the $500$ t-statistics for each pixel, converting them to z-scores to create a null distribution. The blue curve in \Cref{Hist_LST} shows this null distribution.

Next, we compared the Wiki-random and the random datasets by two-sample t-tests on each pixel.
For each pixel, the null hypothesis is that the mean values are the same for the Wiki-random and random datasets.
Our z-scores for these tests are shown in the histogram in \Cref{Hist_LST}. 

\begin{figure}[h!]
    \centering
    \begin{subfigure}[b]{2.5in}
        \centering
        \includegraphics[width=2.3in]{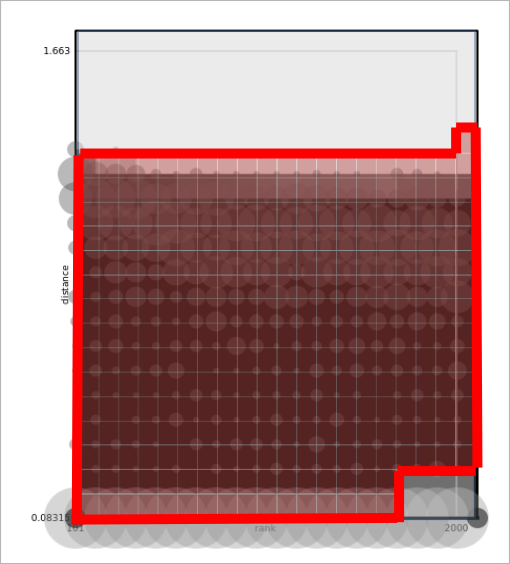}
        \caption{\wikirandom\  dataset}
        \label{rivet:A}
    \end{subfigure}
    \begin{subfigure}[b]{2.5in}
        \centering
        \includegraphics[width=2.3in]{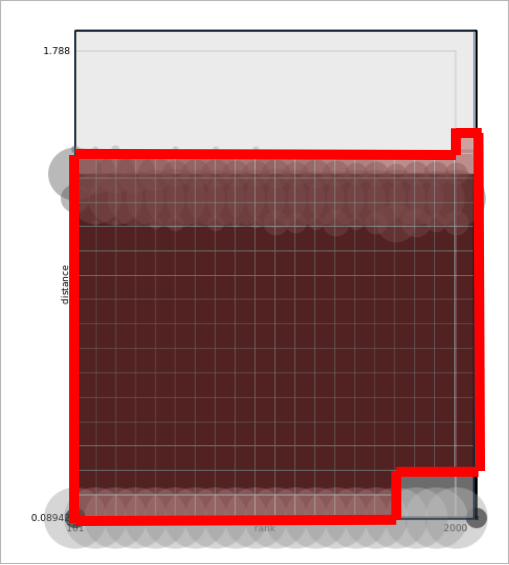}
        \caption{\random\  dataset}
        \label{rivet:B}
    \end{subfigure}
    \caption{Large scale hypothesis testing on $H_0$ homology found 298 significant pixels, which are indicated here by the red outlines on top of the RIVET plots.}
    \label{pixel_pic}
\end{figure}

We found 298 significant pixels out of 400, meaning that their z-scores lie above the 95th percentile from the null distribution.
These significant pixels are shown in the red region in \Cref{pixel_pic}.
This region of the plot shows the range of popularity and scale parameters in which the Wiki-random and random datasets exhibit significantly different numbers of connected components. 
The \emph{power} of a statistical test is the probability of rejecting the null hypothesis when, in fact, it is false. The power of our large scale hypothesis test is the proportion of significant pixels, which is approximately $0.75$. 

Using the same steps as above but with $H_1$ homology, we again compared the Wiki-random and random datasets.
In this large scale hypothesis testing, we 
found  $192$ significant pixels out of $400$, which again reveals significant differences between the two types of datasets, this time in terms of the number of holes in the Rips complexes constructed from the data at particular scales.

\subsubsection{Two-Sample t-Test on Matching Distance}\label{twoSampleMatching}

The second statistical method that we used is the \emph{two-sample t-test on matching distance}. As discussed in \Cref{twoparam}, the matching distance is a topological notion of distance between two datasets.
In this experiment, we used $100$ Wiki-random data sets \wrand{1}, \wrand{2}, ..., \wrand{100} and $100$ random datasets \rand{1}, \rand{2}, ..., \rand{100}.

To create a null distribution, we computed the matching distances between pairs of $H_0$ persistence modules computed from random datasets, since we assume there is no significant topological difference between random datasets.
Since the t-test requires independence, each random data set can be paired with only one other random data set.
Thus, we formed 50 pairs from the datasets \rand{1}, \rand{2}, ..., \rand{100} and computed the matching distance between each pair, resulting in $50$ matching distances.
We used a bootstrapping process to estimate the null distribution of matching distance from our data. Specifically, we sampled $50$ times with replacement from our $50$ matching distances, computing the average from this sample to obtain a bootstrapped estimator. 
We repeated this to compute $1000$ bootstrapped estimators, creating a bootstrapped distribution, which we regard as the null distribution of matching distances. 

\begin{figure}[t]
    \centering
    \includegraphics[scale=0.6]{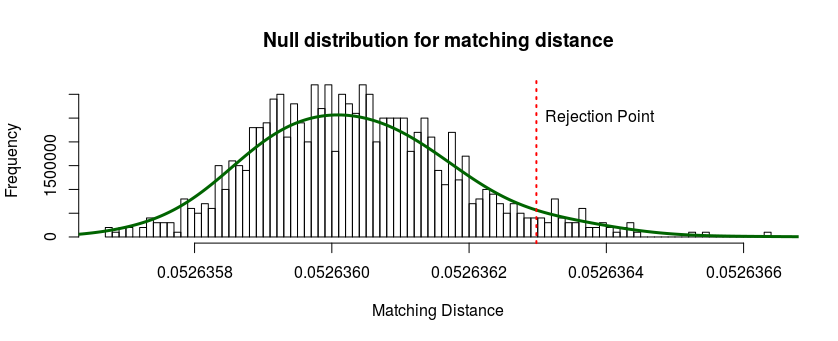}
    \caption{Null distribution of the matching distance is plotted as a histogram and a smooth curve; the red dashed line gives the $95\%$ percentile of the distribution.}
    \label{MD_hist}
\end{figure}

\Cref{MD_hist} shows the histogram of the bootstrapped estimators.
A smoothed version of this histogram is plotted in green, and the red dashed line gives the $95\text{th}$ percentile, which we use as the rejection level for comparing Wiki-random and random data. We note that the null distribution has very small variance about its mean.

We next computed matching distances between Wiki-random and random datasets, again for $H_0$ homology.
We found that the matching distance between 100 pairs of Wiki-random and radom datasets are all $0.368$ to three significant digits, which is far above the rejection level identified earlier and shown in \Cref{MD_hist}.
Thus, we can say that the matching distance between Wiki-random datasets and random datasets differs significantly from the matching distance between pairs of random datasets. We conclude that Wikipedia data sets exhibit different structure from the random datasets, discernible in $H_0$ persistent homology.

In contrast, we were unable to establish statistical significance when repeating this analysis for $H_1$ homology.
The two-sample t-test on matching distance does not reveal statistically significant differences in this case. We suspect this may be due to the coarsening used in our persistence calculation; it would be interesting to re-do the persistence calculation with more bins, as well as the matching distance calculation with more angle and offset values, to see if this affects the results of the t-test.

\subsubsection{Two-Sample t-Test on Bar Lengths}

The third statistical technique that we used is the \emph{two-sample t-test on bar lengths}. 
This is similar to the previous method, but instead of using the matching distance, we examined the barcodes that realize the matching distance.
As explained in \Cref{twoparam}, the matching distance between two persistence modules is a bottleneck distance between two barcodes; we computed the average lengths of bars in this pair of barcodes that realize the matching distance. We regard this average length as a measurement of the average persistence of topological features in the data.

\begin{figure}[t]
    \centering
    \includegraphics[scale=0.6]{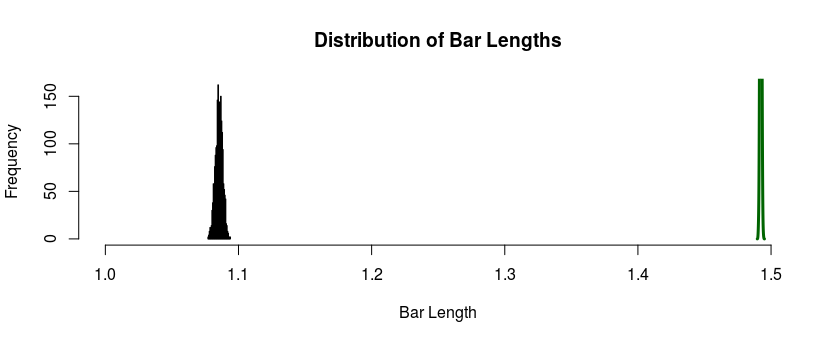}
    \caption{The green curve at right depicts the null distribution for the average $H_0$ bar lengths, and the black histogram shows the average bar lengths for the Wikipedia data.}
    \label{bar_length}
\end{figure}

Similar to the previous method, we computed a null distribution of average bar lengths from the $H_0$ persistent homology of five random datasets, bootstrapping this $1000$ times to compute a null distribution.
The resulting null distribution is approximately normal with a mean of $1.492$ and a standard deviation of $0.0006$. This distribution appears as the green spike in \Cref{bar_length}.

We then computed the barcodes that realize the matching distance between pairs of Wiki-random and random datasets. We found that the average bar lengths from the Wiki-random barcodes have an approximately normal distribution with mean $1.085$ and standard deviation $0.003$, as shown in \Cref{bar_length}. Thus, the average bar length allows us to easily distinguish between the Wiki-random and random datasets.

We also noticed that the $100$ longest bars arising from the random datasets are longer than all the bars arising from the Wiki-random datasets. 
This is surprising given that the points in the random datasets are not farther from the origin, on average, than the points in the Wiki-random datasets.
However, the $100$ longest bars arising from the random datasets tell us that the points in the random datasets are farther apart from each other than are points in the Wiki-random datasets. This is somewhat hard to visualize since these points lie in $\mathbb{R}^{200}$, but we note that high dimensionality makes it possible for a large collection of points to be unit distance apart from each other, while still inside a unit ball. We infer that the random data points are more evenly distributed than the Wiki-random data points.
We observe that we can distinguish between Wiki-random and random data simply by looking at the top bar lengths along the line which realizes matching distance. 

Turning to the $H_1$ persistent homology, we did not find a significant difference between the Wiki-random and random datasets.
The power of this test is the probability of correctly rejecting null hypothesis if we assume that the bar lengths arising from the Wiki-random datasets are longer than the bar lengths arising from the random datasets. 
We found that power of this analysis is $0.322$, which implies that this analysis often fails to distinguish a non-random datasets from our null distribution of average $H_1$ bar lengths.
Also, we found that 20 of the longest bars from the random datasets are longer than the bars from the Wikipedia datasets, which is much fewer than in $H_0$ analysis. 
It appears that the difference between these two types of data is harder to discern when analyzing holes, and easier to discern when analyzing connected components.

\subsection{Distribution of Matching Distance}

Given that the Wikipedia datasets exhibit different topological structure than the random datasets, we now investigate how the topological structure of the Wikipedia datasets change as we replace some of their vectors with vectors from the random datasets. In essence, we replace some of the Wikipedia signal with noise, and we investigate the ability of the matching distance to detect this noise.

We used our \max\ and \random\ datasets to conduct this experiment. 
We first computed the matching distance (using degree-0 homology) between these two datasets. 
Then we selected a vector from \max\ at random and replaced it with a vector from \random, obtaining a new dataset that we call \textsf{MaxRp1}.
We computed the matching distance between \textsf{MaxRp1} and \max. 
We proceeded to replace more Wikipedia vectors with random vectors until only $50$ Wikipedia vectors remained in the dataset. 

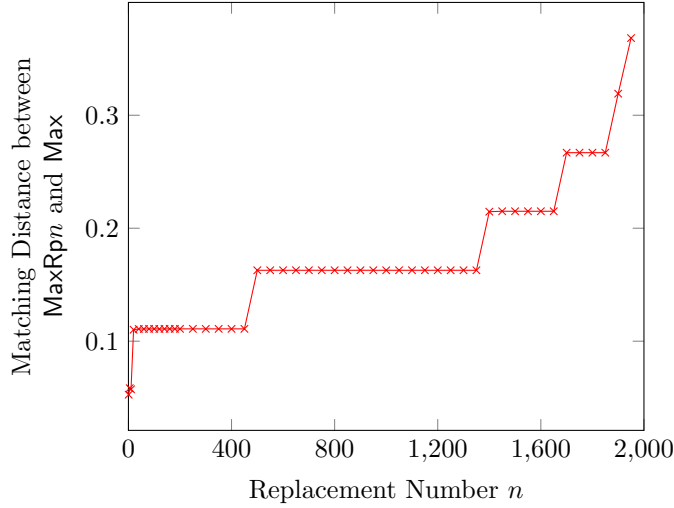
\begin{figure}[ht]
    \centering
\begin{tikzpicture}
\begin{axis}[xmin=0, xmax=2000,xtick = {0,400,800,1200,1600,2000},ytick={0,0.1,0.2,0.3,0.4},xlabel=Replacement Number $n$, ylabel style={align=center},ylabel={Matching Distance between\\\textsf{MaxRp$n$} and \textsf{Max}}]

\addplot[color=red,mark=x] coordinates {(1,0.05263446864)(5,0.05842528485) (10,0.05744390403) (20,0.1101542249) (40,0.1108197802) (60,0.1108443354) (80,0.1109104706) (100,0.1109288877) (120,0.1108443354) (140,0.1108948986) (160,0.1109233195) (180,0.1109071737) (200,0.1109198794) (250,0.1109198794) (300,0.1109119947) (350,0.1109119947) (400,0.1109214004) (450,0.1109132089) (500,0.1627700091) (550,0.1627746893)(600,0.1627746893) (650,0.1627788976) (700,0.1627779542)(750,0.162776907) (800,0.162776907) (850,0.1627667028) (900,0.1627788976) (950,0.1627734314) (1000,0.1627734314) (1050,0.1627746893) (1100,0.1627770193) (1150,0.1627777946) (1200,0.1627777946) (1250,0.1627788976) (1300,0.1627746893) (1350,0.1627788976) (1400,0.2146297285) (1450,0.2150093251) (1500,0.2150093251) (1550,0.2150072151) (1600,0.2150093251) (1650,0.21500767) (1700,0.2668596001) (1750,0.2668597071) (1800,0.2668596001) (1850,0.2668596001) (1900,0.3190906936)(1950,0.3682807685)};

\end{axis}
\end{tikzpicture}
\caption{Plot of matching distance as a function of replacement number when using 20 bins.}
    \label{RP_20}
\end{figure}

\Cref{RP_20} plots the matching distance between \textsf{MaxRp}$n$ and \textsf{Max} as a function of $n$, the number of vectors in \max\ that we replaced with random vectors.
Since the average matching distance among random data sets is about $0.05$ (see \Cref{MD_hist}), we regard a matching distance of $0.05$ as too small to indicate a topological difference between datasets.
We found that a replacement number of less than $20$ resulted in a matching distance of about $0.05$ between the replacement dataset and \max, indicating that we do not detect a topological difference between the underlying point clouds.

\begin{figure}[h!]
    \centering
\begin{tikzpicture}
\begin{axis}[xmin=0, xmax=2000,xtick = {0,400,800,1200,1600,2000},ytick={0,0.1,0.2,0.3,0.4},xlabel=Replacement Number $n$, ylabel style={align=center},ylabel={Matching Distance between\\\textsf{MaxRp$n$} and \textsf{Max}}]

\addplot[color=red,mark=x] coordinates {(1,0.05675972036)(5,0.08117526604) (10,0.08186803306) (20,0.1079988912) (40,0.1078394127) (60,0.1079011702) (80,0.09992893884) (100,0.1078396809) (120,0.1078876223) (140,0.1079135386) (160,0.1079005807) (180,0.1080939454) (200,0.1079092193) (250,0.133541563) (300,0.107903518) (350,0.1234514693) (400,0.1080889513) (450,0.1335416483) (500,0.1331825726) (550,0.1331825726) (600,0.1331865086) (650,0.1331825726) (700,0.1331836602) (750,0.1584511386) (800,0.1584375185) (850,0.1584513555) (900,0.1584473646) (950,0.1584476993) (1000,0.1584473646) (1050,0.1584475501) (1100,0.1584511386) (1150,0.1584511386) (1200,0.1584513555) (1250,0.1837140645) (1300,0.1837172464) (1350,0.1837140542) (1400,0.1838994197) (1450,0.1838994197) (1500,0.2091637404) (1550,0.2091652736) (1600,0.20916519) (1650,0.2344318194) (1700,0.2346169541) (1750,0.2598832401) (1800,0.2598832401) (1850,0.3107844048) (1900,0.3362345207) (1950,0.3362345207)};

\end{axis}
\end{tikzpicture}
\caption{Plot of matching distance as a function of replacement number when using 40 bins.}
    \label{RP_40}
\end{figure} 

As we replace more vectors, the matching distance increases to about $0.11$, where it remains until the replacement number gets to $450$. 
The matching distance again jumps to about $0.21$ for replacement numbers between $450$ and $1400$, and so on, as shown in \Cref{RP_20}.
The matching distance between \max\ and \random\ is $0.368$. However, when we replace $1900$ vectors from \max\  with random vectors, we observe a matching distance less than $0.368$. Thus, we discern that only $100$ Wikipedia vectors out of $2000$ vectors can affect the topological structure of a random data set.

We repeated the experiment using $40$ bins (instead of $20$) in the RIVET calculations (see \Cref{twoparam}). This resulted in a similar plot of matching distances, as shown in \Cref{RP_40}. However, we observe that the jumps between the matching distance values is now about $\frac{1}{40}$. We conjecture that the intervals between observed matching distance values are determined by the number of bins used in the RIVET calculation.

Finally, we compared the replacement datasets with \random. We computed the matching distances between \textsf{maxRp}$n$ and \random. Unsurprisingly, starting with our \max\ dataset and replacing vectors with random vectors causes the matching distance to decrease as the number of replacements increases.
Notably, the matching distance remains about $0.34$ until we introduce $1400$ random vectors. In other words, with $600$ Wikipedia vectors and $1400$ random vectors in the replacement dataset, we get nearly the same matching distance as we compare our \max\ and \random\ datasets. 
This indicates that the presence of $600$ Wikipedia vectors allows us to clearly distinguish the topological structures of replacement dataset from that of \textsf{Random} dataset.

\subsection{Stability of Dimension Plots}

Considering our three statistical methods in both $H_0$ and $H_1$ persistent homology, we observe that large scale hypothesis testing is the most effective method for distinguishing Wikipedia datasets from random datasets. To take our analysis further, we tested the stability of the Hilbert function plots as we replace a Wikipedia vectors with random vectors.

Using our Wiki-random datasets \wrand{1}, ..., \wrand{15}, we computed a null distribution as described in \Cref{largeScaleTesting}. We chose one dataset as the original, and then created the replacement datasets. Following the same steps described in Section 3.3.1, we randomly replaced Wiki-random vectors with random vectors, producing 30 replacement datasets, ranging from one to 30 replacements.
Large-scale hypothesis testing resulted in 234 significant pixels. These pixels are located in the lower-right corner and in a middle horizontal band in the RIVET plot, in contrast with the pattern shown in \Cref{pixel_pic}. We also note that the matching distances between \wikirandom\ and its replacement datasets are relatively stable: the matching distance is about $0.05$ when we replace three or fewer vectors, and the matching distance increases to about $0.1$ when we replace four to thirty vectors.

We performed a similar experiment using the random datasets as the null distribution. This time, we created replacement datasets by replacing random vectors with Wiki-random vectors. We suspected that we would find less stability in this case. Since Wiki vectors contain more information than random vectors, including a few Wiki vectors in an otherwise random dataset should result in different topological structure.
However, large scale testing only revealed $155$ significant pixels, all located in the lower half of the RIVET plot.

We see that starting with a collection of Wiki vectors, adding a few random vectors results in discernible topological change. However, starting with a collection of random vectors and adding a few Wiki vectors does not result in discernible topological change.
We conclude that replacing Wiki vectors with random vectors results in a more significant topological change than replacing random vectors with Wiki vectors.

\section{Discussion and Further Research}

In summary, we examined three statistical measures on two-parameter persistence and demonstrated their applicability using data from Simple English Wikipedia. 
Of the three methods, large scale hypothesis testing is the easiest to compute, since it is based on the Hilbert function of the persistence modules, while the other methods require computing matching distances.
While the matching distance is more computationally intensive, the barcodes involved in the matching distance give deeper insight into the persistence of topological features in the data.
All three techniques are able to distinguish the Wikipedia datasets from random datasets in $H_0$ analysis, but only large scale hypothesis testing on Hilbert function values showed significant differences between the Wikipedia datasets and random datasets in $H_1$ analysis.

We note that our methodology is limited by the amount of coarsening involved in our two-parameter persistence computations.
Specifically, using more bins when computing the two-parameter persistence modules, and using more angle and offset values in the matching distance calculation, would give more accurate results, but at computational cost.

For future work, we would like to explore the robustness of these statistical methods and apply them to other datasets. 
Recent algorithmic advances provide more accurate or faster matching distance computations; it would be interesting to apply these algorithms to data. Specifically, \cite{KLO} shows that the matching distance between two-parameter persistence modules can be computed exactly in polynomial time, and \cite{KerberNigmetov} provides a fast algorithm for approximating the matching distance.
Furthermore, we would like to identify which points in the datasets correspond with the topological features that we determine to be significant. This would tell us, for example, which Wikipedia articles form clusters or cycles that we identify in persistent homology.
Unfortunately, RIVET is not able to identify specific points in the point cloud that contribute to  topological features detected in persistent homology calculations, so we are unable to make such inference at present.

\section*{Acknowledgements}

We thank Shilad Sen at Macalester College for providing the Wikipedia data. This work was supported by NSF DMS-1606967. This work was carried out while Xiaojun Zheng was an undergraduate student at St.\ Olaf College. We are grateful to the Collaborative Undergraduate Research and Inquiry (CURI) program for facilitating undergraduate research at St.\ Olaf College. 
We also thank the anonymous reviewers who gave valuable feedback on a previous version of this paper.

\medskip

\bibliographystyle{plainurl}
{\footnotesize \bibliography{references} }

\begin{thebibliography}{10}

\bibitem{BobrowskiKahle}
Omer Bobrowski and Matthew Kahle.
\newblock Topology of random geometric complexes: a survey.
\newblock {\em Journal of Applied and Computational Topology}, 1(3):331--364,
  June 2018.
\newblock \href {https://doi.org/10.1007/s41468-017-0010-0}
  {\path{doi:10.1007/s41468-017-0010-0}}.

\bibitem{CarlssonZomorodian}
Gunnar Carlsson and Afra Zomorodian.
\newblock The theory of multidimensional persistence.
\newblock {\em Discrete and Computational Geometry}, 42(1):71--93, 2009.
\newblock \href {https://doi.org/10.1007/s00454-009-9176-0}
  {\path{doi:10.1007/s00454-009-9176-0}}.

\bibitem{CSGO}
Frédéric Chazal, Vin Silva, Marc Glisse, and Steve Oudot.
\newblock {\em The Structure and Stability of Persistence Modules}.
\newblock 07 2012.
\newblock \href {https://doi.org/10.1007/978-3-319-42545-0}
  {\path{doi:10.1007/978-3-319-42545-0}}.

\bibitem{Efron}
B.~Efron and T.~Hastie.
\newblock {\em Computer Age Statistical Inference}.
\newblock Cambridge University Press, 2016.
\newblock \href {https://doi.org/10.1017/CBO9781316576533}
  {\path{doi:10.1017/CBO9781316576533}}.

\bibitem{ghristbarcodes}
Robert Ghrist.
\newblock Barcodes: The persistent topology of data.
\newblock {\em Bulletin of the American Mathematical Society}, 45(1):61--75,
  2008.
\newblock \href {https://doi.org/10.1090/S0273-0979-07-01191-3}
  {\path{doi:10.1090/S0273-0979-07-01191-3}}.

\bibitem{HanOkoYadZhe}
So~Mang Han, Taylor Okonek, Nikesh Yadav, and Xiaojun Zheng.
\newblock Distributions of matching distances in topological data analysis.
\newblock {\em SIAM Undergraduate Research Online}, 13, 2020.
\newblock \href {https://doi.org/10.1137/18S017302}
  {\path{doi:10.1137/18S017302}}.

\bibitem{hatcher2002algebraic}
A.~Hatcher.
\newblock {\em Algebraic Topology}.
\newblock Algebraic Topology. Cambridge University Press, 2002.
\newblock \href {https://doi.org/10.1017/S0013091503214620}
  {\path{doi:10.1017/S0013091503214620}}.

\bibitem{hera}
{Hera}: A repository contains software to compute bottleneck and wasserstein
  distances between persistence diagrams.
\newblock \url{https://bitbucket.org/grey_narn/hera}, 2019.

\bibitem{KLO}
Michael Kerber, Michael Lesnick, and Steve Oudot.
\newblock {Exact Computation of the Matching Distance on 2-Parameter
  Persistence Modules}.
\newblock In {\em 35th International Symposium on Computational Geometry (SoCG
  2019)}, volume 129, pages 46:1--46:15, 2019.
\newblock \href {https://doi.org/10.4230/LIPIcs.SoCG.2019.46}
  {\path{doi:10.4230/LIPIcs.SoCG.2019.46}}.

\bibitem{jea_hera}
Michael Kerber, Dmitriy Morozov, and Arnur Nigmetov.
\newblock Geometry helps to compare persistence diagrams.
\newblock {\em Journal of Experimental Algorithmics (JEA)}, 22:1--20, 2017.
\newblock \href {https://doi.org/10.1145/3064175} {\path{doi:10.1145/3064175}}.

\bibitem{KerberNigmetov}
Michael Kerber and Arnur Nigmetov.
\newblock {Efficient Approximation of the Matching Distance for 2-Parameter
  Persistence}.
\newblock In {\em 36th International Symposium on Computational Geometry (SoCG
  2020)}, volume 164, pages 53:1--53:16, 2020.
\newblock \href {https://doi.org/10.4230/LIPIcs.SoCG.2020.53}
  {\path{doi:10.4230/LIPIcs.SoCG.2020.53}}.

\bibitem{Landi2018}
Claudia Landi.
\newblock The rank invariant stability via interleavings.
\newblock In Erin~Wolf Chambers, Brittany~Terese Fasy, and Lori Ziegelmeier,
  editors, {\em Research in Computational Topology}, pages 1--10. Springer,
  2018.
\newblock \href {https://doi.org/10.1007/978-3-319-89593-2_1}
  {\path{doi:10.1007/978-3-319-89593-2_1}}.

\bibitem{lesnick}
Michael Lesnick.
\newblock The theory of the interleaving distance on multidimensional
  persistence modules.
\newblock {\em Foundations of Computational Mathematics}, 15(3):613--650, 2015.
\newblock \href {https://doi.org/10.1007/s10208-015-9255-y}
  {\path{doi:10.1007/s10208-015-9255-y}}.

\bibitem{LesnickWright}
Michael Lesnick and Matthew Wright.
\newblock {Interactive Visualization of 2-D Persistence Modules}.
\newblock {\em ArXiv e-prints}, December 2015.
\newblock URL: \url{https://arxiv.org/pdf/1512.00180.pdf}, \href
  {http://arxiv.org/abs/1512.00180} {\path{arXiv:1512.00180}}.

\bibitem{rivet}
{RIVET}: The rank invariant visualization and exploration tool.
\newblock \url{https://github.com/rivetTDA/rivet}, 2015-2020.

\bibitem{rivetPython}
{Rivet-Python}: Python api for rivet.
\newblock \url{https://github.com/rivetTDA/rivet-python}, 2018.

\bibitem{Shilad}
Shilad Sen, Anja~Beth Swoap, Qisheng Li, Brooke Boatman, Ilse Dippenaar,
  Rebecca Gold, Monica Ngo, Sarah Pujol, Bret Jackson, and Brent Hecht.
\newblock {Cartograph: Unlocking Spatial Visualization Through Semantic
  Enhancement}.
\newblock {\em the 22nd International Conference on Intelligent User
  Interfaces}, 2017.
\newblock \href {https://doi.org/10.1145/3025171.3025233}
  {\path{doi:10.1145/3025171.3025233}}.

\bibitem{Wasserman}
Larry Wasserman.
\newblock {\em All of Statistics: A Concise Course in Statistical Inference}.
\newblock Springer, 2010.

\end{thebibliography}

\end{document}